\documentclass[11pt]{amsart} \usepackage{amssymb,amsmath}
\def\Esimple{E_{\text{simple}}}

\def\F{\mathbb{F}}
\def\Z{\mathbb{Z}}

\begin{document}
\pagestyle{empty}
%\mainmatter

\title{Fast Elliptic Curve Arithmetic\\ and Improved Weil Pairing Evaluation}
%\titlerunning{Fast Elliptic Curve Arithmetic}

\author{Kirsten Eisentr\"ager, Kristin Lauter, and Peter L.~Montgomery}

%\authorrunning{Eisentr\"ager, Lauter, Montgomery}

%\institute{Department of Mathematics,
%           University of California, 
%           Berkeley, CA 94720
%    \email{eisentra@math.berkeley.edu} \\
%\and
%           Microsoft Research, 
%           One Microsoft Way, Redmond, WA 98052 \\
%    \email{klauter@microsoft.com, petmon@microsoft.com}
%}      

\begin{abstract}
  {\small We present an algorithm which speeds scalar multiplication
    on a general elliptic curve by an estimated 3.8\% to 8.5\% over
    the best known general methods when using affine coordinates.
    This is achieved by eliminating a field multiplication when we
    compute $2P + Q$ from given points $P$, $Q$ on the curve.  We give
    applications to simultaneous multiple scalar multiplication and to
    the Elliptic Curve Method of factorization.  We show how this
    improvement together with another idea can speed the computation
    of the Weil and Tate pairings by up to $7.8\%$.}
\end{abstract}

\maketitle


\section{Introduction }
This paper presents an algorithm which can speed scalar multiplication
on a general elliptic curve, by doing some arithmetic differently.
Scalar multiplication on elliptic curves is used by cryptosystems and
signature schemes based on elliptic curves.  Our algorithm saves an
estimated $3.8\%$ to $8.5\%$ of the time to perform a scalar
multiplication on a general elliptic curve, when compared to the
best-known general methods.  This savings is important because the
ratio of security level to computation time and power required by a
system is an important factor when determining whether a system will
be used in a particular context.

Our main achievement eliminates a field multiplication whenever we are
given two points $P$, $Q$ on an elliptic curve and need $2P + Q$ (or
$2P - Q$) but not the intermediate results $2P$ and $P+Q$.  This
sequence of operations occurs many times when, for example,
left-to-right binary scalar multiplication is used with a fixed or
sliding window size.

Some algorithms for simultaneous multiple scalar multiplication
alternate doubling and addition steps, such as when computing $k_1 P_1
+ k_2 P_2 + k_3 P_3$ from given points $P_1$, $P_2$, and $P_3$. Such
algorithms can use our improvement directly.  We give applications of
our technique to the Elliptic Curve Method for factoring and to
speeding the evaluation of the Weil and Tate Pairings.

The paper is organized as follows.  Section~\ref{background} gives
some background on elliptic curves.  Section~\ref{algorithm} gives a
detailed version of our algorithm.  Section~\ref{analysis} estimates
our savings compared to ordinary left-to-right scalar multiplication
with windowing.  Section~\ref{examples} illustrates the improvement
achieved with an example. It also describes applications to
simultaneous multiple scalar multiplication and the Elliptic Curve
Method for factoring.  Section~\ref{SecPairings} adapts our technique
to the Weil and Tate pairing algorithms. Appendix~\ref{code} gives the
pseudocode for implementing the improvement, including abnormal cases.

\section{Background} \label{background} 
Elliptic curves are used for several kinds of cryptosystems, including
key exchange protocols and digital signature algorithms \cite{IEEE}.
If $q$ is a prime or prime power, we let $\F_q$ denote the field with
$q$ elements.  When $\gcd(q,\, 6) = 1$, an elliptic curve over the
field $\F_{q}$ is given by an equation of the form
$$
\Esimple:\;y^2 = x^3 + ax + b
$$ 
with $a,~b$ in $\F_{q}$ and $4a^3 + 27b^2 \neq 0$. (See
\cite[p.~48]{Silverman}.)

A more general curve equation, valid over a field of any characteristic, 
is considered in Appendix~\ref{code}.  The general curve equation subsumes the 
case
$$
E_2: \;y^2 + xy = x^3 + ax^2 + b
$$
with $a,~b$ in $\F_{q}$ and $b \neq 0$, 
which is used over fields of characteristic~2.
 
In all cases the group used when implementing the cryptosystem is the
group of points on the curve over $\F_q$.  If represented in affine
coordinates, the points have the form: $(x,~y)$, where $x$ and $y$ are
in $\F_q$ and they satisfy the equation of the curve, as well as a
distinguished point~$\mathbf{O}$ (called the {\it point at infinity})
which acts as the identity for the group law.  Throughout this paper we
work with affine coordinates for the points on the curve.
 
Points are added using a geometric group law which can be expressed
algebraically through rational functions involving $x$ and $y$.
Whenever two points are added, forming $P+Q$, or a point is doubled,
forming $2P$, these formulae are evaluated at the cost of some number
of multiplications, squarings, and divisions in the field.  For
example, using $\Esimple$, to double a point in affine
coordinates costs $1$~multiplication, $2$~squarings, and $1$~division
in the field, not counting multiplication by $2$ or~$3$
\cite[p.~58]{BSS}.  To add two {\it distinct} points in affine
coordinates costs $1$~multiplication, $1$~squaring, and $1$~division
in the field.  Performing a doubling and an addition $2P+Q$ costs~$2$
multiplications, $3$~squarings and $2$~divisions if the points are
added as $(P+P)+Q$, i.e., first double $P$ and then add $Q$.

\section{The Algorithm} \label{algorithm}
Our algorithm performs a doubling and an addition, $2P+Q$, on an
elliptic curve $\Esimple$ using only $1$~multiplication,
$2$~squarings, and $2$~divisions (plus an extra squaring when $P =
Q$).  This is achieved as follows: to form $2P+Q$, where $P =
(x_1,~y_1)$ and $Q = (x_2,~y_2)$, we first find $(P + Q)$, except we
omit its $y$-coordinate, because we will not need that for the next
stage.  This saves a field multiplication.  Next we form $(P + Q) +
P$.  So we have done two point additions and saved one
multiplication. This trick also works when $P = Q$, i.e., when tripling
a point.  One additional squaring is saved when $P \neq Q$ because
then the order of our operations avoids a point doubling.

Elliptic curve cryptosystems require multiplying a point $P$
by a large number~$k$.  If we write $k$ in binary form and compute
$k P$ using the left-to-right method of binary scalar multiplication, we can
apply our trick at each stage of the partial computations.

Efficient algorithms for group scalar multiplication have a long history (see
\cite{Knuth} and \cite{Gordon}), and optimal scalar multiplication routines
typically use a combination of the left-to-right or right-to-left
$m$-ary methods with sliding windows,
addition-subtraction chains, signed representations, etc.  Our
procedure can be used on top of these methods for $m=2$ to obtain a
savings of up to $8.5\%$ of the total cost of the scalar multiplication for
curves over large prime fields, depending upon the window size and form which
is used.  This is described in detail in Section~\ref{analysis}.

%This can be seen from the fact that a
%field inversion costs roughly the equivalent of 4 field
%multiplications for large prime fields of size of cryptographic
%interest (such as 163-bit fields), and a field division does one
%inversion and one multiplication.  Throughout this paper, we also
%estimate that a field squaring costs roughly the same as a field
%multiplication, which is true for such odd-characteristic fields.
%Therefore, saving one field multiplication accounts for about $6\%$ of
%the total cost of exponentiation, plus we save a squaring due to the
%order of operations.  The examples in the next section give some
%relatively small exponentiations to illustrate how this works.

\subsection{Detailed Description of the Algorithm}

Here are the detailed formulae for our procedure when the curve has the
form $\Esimple$ and all the points are distinct, none equal
to $\mathbf{O}$.  Appendix~\ref{code} gives details for all
characteristics.  That appendix also covers special cases, where an
input or an intermediate result is the point at infinity.

Suppose $P=(x_1,~y_1)$ and $Q=(x_2,~y_2)$ are
distinct points on $\Esimple$, and $x_1 \neq x_2$.
The point $P+Q$ will have coordinates $(x_3,~y_3)$, where
\begin{eqnarray*}
\lambda_1 &=& (y_2-y_1)/(x_2-x_1), \\ x_3 &=& \lambda_1^2 - x_1 -
x_2, \qquad \mbox{ and }\\ y_3 &=& (x_1-x_3)\lambda_1 - y_1.
\end{eqnarray*}

Now suppose we want to add $(P+Q)$ to $P$.  
We must add $(x_1,~y_1)$ to $(x_3,~y_3)$ using the above rule.  
Assume $x_3 \neq x_1$. 
The result has coordinates $(x_4, ~y_4)$, where
\begin{eqnarray*}
\lambda_2  &=& (y_3-y_1)/(x_3-x_1), \\
x_4 &=& \lambda_2^2 -x_1 -x_3, \qquad \mbox{ and }\\
y_4 &=& (x_1 -x_4)\lambda_2  - y_1.
\end{eqnarray*}

We can omit the $y_3$ computation, because it is used only
in the computation of $\lambda_2$, which can be computed
without knowing $y_3$ as follows:
$$
\lambda_2  = -\lambda_1 - 2y_1/(x_3-x_1).
$$
Omitting the $y_3$ computation saves a field multiplication.
Each $\lambda_2$ formula requires a field division, so the
overall saving is this field multiplication.  

This trick can also be applied to save one multiplication when
computing $3P$, the triple of a point $P \neq \mathbf{O}$, where the
$\lambda_2$ computation will need the slope of a line through two
distinct points $2P$ and $P$.

This trick can be used twice to save $2$~multiplications when
computing $3P + Q = ((P + Q) + P) + P$. 
Thus $3P + Q$ can be computed using
$1$~multiplication, $3$~squarings, and~$3$ divisions.
Such a sequence of operations would be performed repeatedly 
if a multiplier
were written in ternary form and left-to-right scalar multiplication were
used.  Ternary representation performs worse than binary
representation for large random multipliers~$k$, but the operation
of triple and add might be useful in another context.

A similar trick works for elliptic curve arithmetic
in characteristic~2, as
is shown in the pseudocode in Appendix~\ref{code}.

Table~\ref{tabcosts} summarizes the costs of some operations on
$\Esimple$.

\begin{table}[h]
\caption{Costs of simple operations on $\Esimple$}\label{tabcosts}
\begin{center}
\renewcommand{\arraystretch}{1.1}
\begin{tabular}{|c|c|l|}      
\hline
       Doubling   & $2P$    &  2 squarings, 1 multiplication, 1 division  \\
%\hline
       Add        &  $P \pm Q$  &  1 squaring,\phantom{s}  1 multiplication, 1 division  \\
%\hline
       Double-add & $2P\pm Q$  &  2 squarings, 1 multiplication, 2 divisions \\
%\hline
       Tripling   & $3P$       &  3 squarings, 1 multiplication, 2 divisions \\
%\hline
       Triple-add & $3P\pm Q$  &  3 squarings, 1 multiplication, 3 divisions \\
\hline
\end{tabular}
\end{center}

\end{table}

\section{Comparison to Conventional Scalar Multiplication} \label{analysis}
In this section we analyze the performance of our algorithm compared
to conventional left-to-right scalar multiplication. We will refer to adding
two distinct points on the curve $E$ as elliptic curve addition, and
to adding a point to itself as elliptic curve doubling. 
Suppose we would like to compute $k P_0$ given $k$ and
$P_0$, where the exponent~$k$ has $n$~bits and $n$ is at least $160$.

%ratio

Assume that the relative costs of field operations are $1$~unit per
squaring or general multiplication and $\alpha$~units per inversion.
\cite[p.~72]{BSS} assumes that the cost of an inversion is between $3$
and $10$ multiplications.  In some implementations the relative cost
of an inversion depends on the size of the underlying field.  Our own
timings on a Pentium~II give a ratio of 3.8 for a $160$-bit prime
field and 4.8 for a $256$-bit prime field when not using Montgomery
multiplication.  Some hardware implementations for fast execution of
inversion in binary fields yield inversion/multiplication ratios of
$4.18$ for $160$-bit exponents and $6.23$ for $256$-bit exponents
\cite{Koc2002}.

%%%%%{\tt \bf PLM -- Is `exponents' the proper word in the Koc reference?
%  Or do you want the degree of the field?}
% KL I think exponents is the right choice, since they refer to bitsize
% in their table.
%%KE Kristin, should we mention the hardware issue in the Koc paper?
%epsilon

The straightforward left-to-right binary method
needs about $n$ elliptic curve doublings. If the window size is one,
then for every $1$-bit in the binary representation, we
perform an elliptic curve doubling followed directly by an elliptic
curve addition.  Suppose about half of the bits in the binary
representation of $k$ are $1$'s. Then forming $k P$ consists of
performing $n$ elliptic curve doublings and $n/2$ elliptic curve
additions.

In general, independent of the window size, the number of elliptic
curve doublings to be performed will be about $n$ asymptotically,
whereas the number of elliptic curve additions to be performed will
depend on the window size.  Define the value $0<\varepsilon<1$ for a
given window size to be such that the number of elliptic curve
additions to be performed is $\varepsilon n$ on average.
For example with window size~1, $\varepsilon$ is $1/2$.

If we fix a window size and its corresponding
$\varepsilon$, then the conventional algorithm 
for scalar multiplication
needs about $2n + \varepsilon n$ field squarings, $n + \varepsilon n$
field general multiplications, and $n + \varepsilon n$ field
divisions. If one inversion costs $\alpha$ multiplications, then the
cost of a division is ($\alpha+1$) multiplications. So the overall
cost in field multiplications is
$$
(2n + \varepsilon n) + (n + \varepsilon n) + (\alpha+1)(n + \varepsilon n)
= (4 + \alpha)n + (3+ \alpha) \varepsilon n .
$$

Now we analyze the percentage savings obtained by our algorithm,
not including precomputation costs. The
above computation includes $\varepsilon n$ sub-computations of the
form $2 P_1 + P_2$.  Writing each as $P_1 + (P_1 + P_2)$ saves one
squaring per sub-computation, reducing the overall cost to $(4 +
\alpha) n + (2 + \alpha) \varepsilon n$.  The technique in 
Section~\ref{algorithm} saves another multiplication per sub-computation,
dropping the overall cost to $(4 + \alpha) n + (1+\alpha) \varepsilon
n$. This means we get a savings of 
$$2\varepsilon/((4 + \alpha)+(3 + \alpha)\varepsilon).$$

When the window size is $1$ and the inversion/multiplication ratio
$\alpha$ is assumed to be $4.18$, this gives a savings of 
8.5\%. When $\alpha$ is assumed to be $6.23$, we still obtain a savings of 
$6.7\%$.
When the window size is $2$ and  $2P$ and $3P$ have been precomputed, we find
that $\varepsilon=3/8$.  So when $\alpha$ is $4.18$, we get a savings
of 6.9\%, and when $\alpha$ is $6.23$, we still obtain a savings of
$5.5\%$.
Similarly if the window size is $4$, and we have precomputed small multiples
of $P$, we still achieve a savings of $3.8\%$ to $4.8\%$, depending on
$\alpha$.

Another possibility is using addition/subtraction chains and 
higher-radix methods. 
The binary method described in \cite[section
A.10.3]{IEEE} utilizes addition/subtraction chains and does about
$2n/3$ doublings and $n/3$ double-adds (or double-subtracts), so
$\varepsilon = 1/3$ in this case. (See \cite[section 2.3]{Gordon} for
an explanation of how we obtain $\varepsilon = 1/3$ in this case.)
With $\alpha = 4.18$, we get a $6.3\%$ improvement. 

Scalar multiplication algorithms that use addition/sub\-trac\-tion
chains as well as sliding window size may have lower $\varepsilon$, but we
still obtain at least a $4.2\%$ savings if $\varepsilon > 0.2$ and $\alpha
= 4.18$.

\cite[Section 3.3]{SaSa2001} presents some possible trade-offs
arising from different 
inversion/mul\-tipli\-cation ratios.  We discuss this further in 
Section~\ref{ssECM}.

\section{Examples and Applications} \label{examples}

\subsection{Left-to-Right Binary Scalar Multiplication}

Suppose we would like to compute $1133044P = (100010100100111110100)_2
P$ with left-to-right binary method.
% and a fixed window size of $4$. 
We will do this twice, the standard way and the new way.
%The standard way computes $16P+Q$ using $4$
%doublings and $1$ addition; our improved method computes it as
%$(8P+Q)+8P$ with $3$ point doublings and $2$ additions and takes
%advantage of the savings described in the previous section.
For each method, we assume that $3P$ has been precomputed. 
The next table compares the number of operations needed ($a =$ point
additions, $d=$ point doublings, $div=$ field divisions,
$s =$ field squarings, $m =$ field multiplies):
\begin{alignat*}{3}         
   &&{\mbox{Standard }}\qquad &{\mbox{Improved}}\\
% &&\mbox{implementation }\;\;&&\\
  1133044P& = 4(283261P)&2d\qquad &2d\\ 283261P &= 128(2213P) - 3P
  &7d+1a\qquad &6d+2a (\mbox{save } 1m)\\ 2213P &= 8(277P) - 3P
  &3d+1a\qquad &2d+2a (\mbox{save } 1m)\\ 277P &= 8(35P) - 3P
  &3d+1a\qquad &2d+2a (\mbox{save } 1m)\\35P &= 8(4P) + 3P
  &3d+1a\qquad &2d+2a (\mbox{save } 1m)\\4P &= P + 3P &1a\qquad &1a \\
\text{ Total:}& &23 div+41s+23m \qquad &23 div+37s+19m
\end{alignat*}

%\begin{alignat*}{3}         
%  & &{\mbox{standard implementation}}\qquad &{\mbox{improvements}}\\
%  Tmp &:= P\\ Tmp &:= 16Tmp + P &4d+a\qquad &3d+2a (\mbox{save } 1m)\\
%  Tmp &:= 16Tmp + 4P &4d+a\qquad &3d+2a (\mbox{save } 1m)\\ Tmp &:=
%  16Tmp + 9P &4d+a\qquad &3d+2a (\mbox{save } 1m)\\ Tmp &:= 16Tmp + 15P
%  &4d+a\qquad &3d+2a (\mbox{save } 1m)\\ Tmp &:= 16Tmp + 4P &4d+a\qquad
%  &3d+2a (\mbox{save }1m)\\ Total & &25div+45s+25m \qquad
%  &25div+40s+20m
%\end{alignat*}
This saves 4 squarings and 4 multiplications.  Estimating the division
cost at about 5 multiplications, this savings translates to about $4.47\%$. 
% Sliding window sizes can also be used in combination with our
%improvement.

\subsection{Simultaneous Multiple Scalar Multiplication}

Another use of our elliptic curve double-add technique is
multiple scalar multiplication, such as $k_1P_1 + k_2P_2 + k_3P_3$, where the
multipliers $k_1$, $k_2$, and $k_3$ have approximately the same
length.  One algorithm creates an $8$-entry table with
$$\mathbf{O},\quad P_1,\quad P_2,\quad P_2+P_1,\quad P_3, \quad P_3 +
P_1,\quad P_3 +P_2,\quad P_3 +P_2+P_1.
$$
Subsequently it uses one elliptic curve doubling followed by the
addition of a table entry, for each multiplier bit
\cite{Moller2001}.  About $7/8$ of the doublings are followed by
an addition other than $\mathbf{O}$.

    To form $29P_1 + 44P_2$, for example, write the
multipliers in binary form: $(011101)_2$ and $(101100)_2$.
Scanning these left-to-right, the steps are
%\\[0,1cm]
%\begin{table}[h]
\begin{center}
\renewcommand{\arraystretch}{1.1}
\begin{tabular}{|c|c|l|}
\hline
   Bits   &  Table entry   &     Action \\
\hline
      $0, 1$    &    $ P_2$     &         $ T := P_2$  \\
\hline
      $1, 0$    &    $ P_1$       &        $T:= 2T +P_1 = P_1 + 2P_2$  \\
\hline
      $1, 1$    &  $ P_1+P_2$     &     $T := 2T + (P_1+P_2) = 3P_1 + 5P_2$ \\
\hline
      $1, 1$   &   $ P_1+P_2$     &     $T := 2T + (P_1+P_2) = 7P_1 + 11P_2$ \\
\hline
      $0, 0$   &  $\mathbf{O}$    &         $T := 2T = 14P_1 + 22P_2$  \\
\hline
      $1, 0$   &     $ P_1$      &       $T := 2T +P_1 = 29P_1 + 44P_2$  \\
\hline
\end{tabular}
\end{center}
%\end{table}
%
{\noindent There is one elliptic curve addition $(P_1+P_2)$ to construct the
four-entry table, four doublings immediately followed by an addition, and one
doubling without an addition.  While doing $10$~elliptic curve
operations, our technique is used four times.  Doing the
multipliers separately, say by the addition-subtraction chains}
$$ 1, 2, 4, 8, 7, 14, 28, 29 \qquad \text{and} \qquad
1, 2, 4, 6, 12, 24, 48, 44 $$
takes seven elliptic
curve operations per chain, plus a final add ($15$ total).

\subsection{Elliptic Curve Method of Factorization} \label{ssECM}
     The Elliptic Curve Method (ECM) of factoring a composite
integer~$N$ chooses an elliptic curve~$E$ with coefficients
modulo~$N$.  ECM multiplies an initial point~$P_0$ on $E$ by a large integer~$k$, 
working in the ring $\Z/N\Z$ rather than over a field.
ECM may encounter a zero divisor while trying to invert a 
nonzero integer, but that is good, because it leads to a factorization of~$N$.  
ECM uses only the $x$-coordinate of $k P_0$.

\cite[pp.~260ff]{Mont1987} proposes 
a para\-me\-teri\-za\-tion, 
$By^2=x^3 + Ax^2 +x$, which uses no inversions during
a scalar multiplication and omits the $y$-coord\-inate of the result.  
Its associated costs for computing the $x$-coordinate are
% \begin{table}[h]
\begin{center}\renewcommand{\arraystretch}{1.25}
\begin{tabular}{|l|l|}      
\hline
      $P+Q$ from $P$, $Q$, $P-Q$ &    2 squarings, 4 multiplications \\
%\hline
      $2P$ from $P$              &    2 squarings, 3 multiplications \\
\hline
\end{tabular}
\end{center}
% \end{table}
\noindent{To form $kP$ from $P$ for 
a large $n$-bit integer~$k$, this method uses
about $4n$ squarings and $7n$ multiplications,
working from the binary representation of~$k$.
Some variations \cite{MontLucas} use fewer steps
but are harder to program.}

In contrast, using our technique and the
method in \cite[section A.10.3]{IEEE}, we do about 
$2n/3$ doublings and $n/3$ double-adds (or double-subtracts).
By Table~\ref{tabcosts}, the estimated cost of $k P$
is $2n$~squarings, $n$~multiplications and $4n/3$~divisions.

     The new technique is superior if $4n/3$ divisions cost
less than $2n$ squarings and $6n$ multiplications.  A division can be
implemented as an inversion plus a multiplication, so the new
technique is superior if an inversion is cheaper than 1.5 squarings
and 3.5 multiplications.

    \cite{Mont1987} observes that one may trade two independent
inversions for one inversion and three multiplications, using $x^{-1}
= y(xy)^{-1}$ and $y^{-1} = (xy)^{-1}x$.  When using many curves to
(simultaneously) tackle the same composite integer, the asymptotic
cost per inversion drops to 3~multiplications.

\section{Application to Weil and Tate Pairings}\label{SecPairings}
The Weil and Tate pairings are becoming important for public-key
cryptography \cite{Joux2002}.  The algorithms for these pairings
construct rational functions with a prescribed pattern of poles and
zeroes.  An appendix to \cite{BoFr2001} describes Miller's algorithm
for computing the Weil pairing on an elliptic curve in detail.

Fix an integer~$m > 0$ and an $m$-torsion point~$P$ on an elliptic
curve~$E$. Let $f_1$ be any nonzero field element. For an integer $c >
1$, let $f_{c}$ be a function on $E$ with a $c$-fold zero at~$P$, a
simple pole at~$cP$, a pole of order~$c-1$ at~$\mathbf{O}$, and no
other zeroes or poles. When $c=m$, this means that $f_m$ has an
$m$-fold zero at~$P$ and a pole of order~$m$ at~$\mathbf{O}$.
Corollary 3.5 on page 67 of \cite{Silverman} asserts that such a
function exists.  This $f_c$ is unique up to a nonzero multiplicative
scalar.
%much as an antiderivative is unique up to an additive constant. 
Although $f_c$ depends on $P$, we omit the extra subscript~$P$.

The Tate pairing evaluates a quotient of the form
$f_{m}(Q_1)/f_m(Q_2)$ for two points~$Q_1,\,Q_2$ on $E$ (see, for
example, \cite{BKLS2002}).  (The Weil pairing has four such
computations.)  Such evaluations can be done iteratively using an
addition/subtraction chain for $m$, once we know how to construct
$f_{b+c}$ and $f_{b-c}$ from $(f_b,~bP)$ and $(f_c,~cP)$.  Let
$g_{b,c}$ be the line passing through the points $bP$ and $cP$.  When
$bP = cP$, this is the tangent line to $E$ at $bP$.  Let $g_{b+c}$ be
the vertical line through $(b+c)P$ and $-(b+c)P$.  Then we have the
useful formulae
$$
    f_{b+c}=f_b \cdot f_c \cdot \frac{g_{b,c}}{g_{b+c}} \qquad
\text{and} \qquad f_{b-c}=\frac{f_b \cdot g_{b} }{f_c \cdot g_{-b,c}}.
$$ 

Denote $h_b = f_b(Q_1) / f_b(Q_2)$ for each integer $b$. 
Although $f_b$ was defined only up to a multiplicative 
constant, $h_b$ is well-defined.  We have
\begin{equation} \label{hbc}
    h_{b+c}= h_b \cdot h_c \cdot 
               \frac{g_{b,c}(Q_1) \cdot g_{b+c}(Q_2)}
                        {g_{b,c}(Q_2) \cdot g_{b+c}(Q_1)}
\quad \text{and} \quad
    h_{b-c}=\frac{h_b \cdot g_{b}(Q_1) \cdot g_{-b,c}(Q_2) }
                     {h_c \cdot g_{b}(Q_2) \cdot g_{-b,c}(Q_1) }.
\end{equation}
So far in the literature, only the $f_{b+c}$ formula appears, 
but the $f_{b-c}$ formula is  useful if using addition/subtraction
chains.  
%It is important to note that the multiples of the point $P$,
%$bP$ and $cP$, are assumed to be known inputs to this computation.
The addition/subtraction chain iteratively 
builds $h_m$ along with $mP$.

\subsection{Using the Double-Add Trick with Parabolas}

We now describe an improved method for obtaining $(h_{2b+c},~(2b+c)P)$
given $(h_b,~bP)$ and $(h_c,~cP)$.  The version of Miller's algorithm
described in \cite{BKLS2002} uses a left-to-right binary method
with window size one.  That method would first compute $(h_{2b},
~2bP)$ and later $(h_{2b+c}, ~(2b+c)P)$.  We propose to compute
$(h_{2b+c}, ~(2b+c)P)$ directly, producing only the $x$-coordinate of
the intermediate point $bP+cP$.  To combine the two steps, we
construct a parabola through the points $bP$, $bP$, $cP$, $-2bP-cP$.

To form $f_{2b+c}$, we form $f_{b+c}$ and $f_{b+c+b}$.
The latter can be expressed as 
$$
f_{2b+c}
= f_{b+c} \cdot \frac{f_b \cdot g_{b+c,b}}{g_{2b+c}} 
=  \frac{ f_b \cdot f_c \cdot g_{b,c} } {g_{b+c}}
       \cdot
            \frac{ f_b \cdot g_{b+c,b}} {g_{2b+c}}
     =  \frac{ f_b \cdot f_c \cdot f_b} {g_{2b+c}}
          \cdot 
     \frac{g_{b,c}\cdot g_{b+c,b}}{g_{b+c}}.
$$
We replace $(g_{b,c}\cdot g_{b+c,b})/g_{b+c}$ by the parabola, whose
formula is given below.  Evaluate the formula for $f_{2b+c}$ at $Q_1$
and $Q_2$ to get a formula for $h_{2b+c}$.

\subsection{Equation for Parabola Through Points}
If $R$ and $S$ are points on an elliptic curve~$E$,
then there is a (possibly degenerate) parabolic equation passing
through $R$ twice (i.e., tangent at $R$) 
and also passing through $S$ and $-2R - S$.  
Using the notations $R = (x_1,~y_1)$ and $S =
(x_2,~y_2)$ with $R+S = (x_3, ~y_3)$ and $2R + S = (x_4, y_4)$, 
a formula for this parabola is
\begin{equation}\label{parab1}
   \frac{(y + y_3 - \lambda_1 (x - x_3) ) (y - y_3 - \lambda_2 (x - x_3) )}
        { x - x_3}. 
\end{equation}
The left half of the numerator of (\ref{parab1}) is a line passing
 through $R$, $S$, and $-R-S$ whose slope is $\lambda_1$.  The right
 half of the numerator is a line passing through $R+S$, $R$, and $-2R
 - S$, whose slope is $\lambda_2$.  The denominator is a (vertical)
 line through $R+S$ and $-R-S$.  The quotient has zeros at $R$, $R$,
 $S$, $-2R-S$ and a pole of order four at $\mathbf{O}$.

We simplify (\ref{parab1}) by expanding it in powers of $x - x_3$.  
Use the equation for $E$ 
to eliminate references to $y^2$ and $y_3^2$.     
\begin{equation}
\begin{array}{rl}
    {} & \phantom{=} \frac{y^2 - y_3^2}{x - x_3} - \lambda_1 (y - y_3)
               - \lambda_2(y + y_3) + \lambda_1 \lambda_2 (x - x_ 3)
               \\ {} &= x^2 + x x_3 + x_3^2 + a + \lambda_1 \lambda_2
               (x - x_3) - \lambda_1 (y - y_3) - \lambda_2 (y + y_3)
               \\ {} &= x^2 + (x_3 + \lambda_1 \lambda_2)x -
               (\lambda_1 + \lambda_2) y +
               \text{constant}.\label{parab2}
\end{array}
\end{equation}    
Knowing that (\ref{parab2}) passes through $R = (x_1, y_1)$, one
formula for the parabola is
\begin{eqnarray}\label{parabola}
(x - x_1) (x + x_1 + x_3 + \lambda_1 \lambda_2) 
- (\lambda_1 + \lambda_2) (y - y_1) .
\end{eqnarray}
In the previous section we can now replace $(g_{b,c}\cdot
g_{b+c,b})/g_{b+c}$ by the parabola (\ref{parabola}) with $R=bP$ and
$S=cP$.

 Formula~(\ref{parabola}) for the parabola does not reference $y_3$
and is never identically zero since its $x^2$ coefficient is $1$.
Appendix~\ref{code} gives a formula for this parabola in degenerate cases, as
well as for a more general curve.

\subsection{Savings}

We claim the pairing algorithm needs less effort to evaluate a
parabola at a point than to evaluate lines and take their product at
that point.  The parabola does not reference $y_3$, so we can omit
the $y$-coordinate of $bP+cP$ and can use the double-add trick.

Here is a precise analysis of the savings we obtain by using
the parabola when computing the Tate pairing. 
Again assume that we use the binary method in \cite[section
A.10.3]{IEEE} to form $mP$, where $m$ has $n$ bits.
(It does $2n/3$ doublings and
$n/3$ double-adds or double-subtracts.)
We manipulate the numerator and denominator 
of $h_j$ separately, doing one division
$h_j= h_{\text{num},j} / h_{\text{denom},j}$ at the very end. 

{\bf  Analysis of doubling step:} The analysis of the doubling step is the
same in the standard and in the new algorithms. 
Suppose we want to
compute $(h_{2b},~2bP)$ from $(h_b,~bP)$. We need an elliptic
curve doubling to compute $2(bP)$, after which we apply~(\ref{hbc}).  
If $bP = (x_1,\,y_1)$ and $2bP = (x_4,\,y_4)$ then 
\begin{equation}\label{gbbg2b}
\frac{g_{b,b}}{g_{2b}} 
   = \frac{y - y_1 - \lambda_1(x - x_1)}
          {x - x_4} .
\end{equation}
The doubling (including $\lambda_1$ computation)
costs 3 multiplications and a division.
Evaluating (\ref{gbbg2b}) at $Q_1$ and $Q_2$
(as fractions) costs 2 multiplications.  
Multiplying four fractions in (\ref{hbc}) costs $6$ multiplications.
The net cost is $3 + 2 + 6 = 11$ field multiplications 
(or squarings) and a field division.

{\bf Analysis of double-add step:} 
The standard algorithm performs
one doubling followed by an addition 
to compute $(h_{2b+c},~(2b+c)P)$
from $(h_{b},~bP)$ and $(h_{c},~cP)$. 
Similar to the above analysis we can
compute the cost as $21$ field multiplications and $2$ divisions.
[The cost would be one fewer multiplication if one does 
two elliptic curve additions: $(2b+c)P = (bP + cP) + bP$.]

The new algorithm does one elliptic curve double-add operation.  It
costs only one multiplication to construct the coefficients of the
parabola~(\ref{parabola}), because we computed $\lambda_1$ and
$\lambda_2$ while forming $(2b+c)P$.  Evaluating the parabola (and the
vertical line $g_{2b+c}$) twice costs four multiplications.
Multiplying five fractions costs another $8$ multiplications.  The
total cost is $3 + 1 + 4 + 8 = 16$ field multiplications and $2$ field
divisions.

{\bf Total savings:}
Estimating a division as $5.18$ multiplications, the standard
algorithm for $(h_m,~mP)$ takes $(16.18 \cdot 2n/3) + 
(31.36 \cdot n/3)=(21.24)n$ steps, 
compared to $(16.18 \cdot 2n/3)+(26.36 \cdot n/3)=19.57n$ 
steps for the new method, a $7.8\%$ improvement.  
A Weil pairing algorithm using the parabola
will also save $7.8\%$ over Miller's algorithm, 
because we can view the Weil pairing as ``two
applications of the Tate pairing'', each saving $7.8\%$.

    Sometimes (e.g., \cite{BLS}) one does multiple 
Tate pairings 
with $P$ fixed but varying $Q_1$ and $Q_2$.  
If one has precomputed 
all coefficients of the lines and parabolas, 
then the costs of evaluation are 8 multiplications per doubling 
step or addition step, and 12 multiplications per combined 
double-add step.  The overall costs are $32n/3$ 
multiplications per evaluation with the traditional method 
and $28n/3$ multiplications with the parabolas, a $12.5\%$
improvement.

\appendix

\section{Pseudocode} \label{code}

The general Weierstrass form for the equation of an elliptic curve is:
\begin{equation}
     E: y^2 + a_1 xy  + a_3 y = x^3 + a_2 x^2 + a_4 x + a_6 ,\label{WGEN}
\end{equation}
subject to the condition that the coefficients $a_1$, $a_2$, $a_3$,
$a_4$, $a_6$ satisfy a certain inequality to prevent singularity
\cite[p.~46]{Silverman}.  
The negative of a point $P = (x_1, \, y_1)$ on (\ref{WGEN}) is
$-P = (x_1, ~-a_1 x_3 - a_3 - y_1)$.  [This seems to require 
a multiplication $a_1 x_3$, but in practice $a_1$ is 0 or 1.] 
If $P = (x_1, \, y_1)$ is a finite point on (\ref{WGEN}), then the
tangent line at $P$ has slope
\begin{equation} \lambda_1 = 
\frac{3 x_1^2 + 2 a_2 x_1 + a_4 - a_1 y_1}
     {2 y_1 + a_1x_1 + a_3} .\label{TAN}
\end{equation}

Figure~1 gives the pseudocode for implementing the savings for an 
elliptic curve of this general form.  
Given two points $P = (x_1, \, y_1)$ and $Q = (x_2, \, y_2)$
on $E$, it describes how to compute $2P + Q$ as well as the equation for a 
(possibly degenerate) parabola through $P$, $P$, $Q$, and $-(2P+Q)$.

Often the curve coefficients in (\ref{WGEN}) are chosen
to simplify (\ref{TAN}) --- the precise choices depend on the field.
For example, it is common in characteristic~2 \cite[p.~115]{IEEE} to
choose $a_1 = 1$ and $a_3 = a_4 = 0$, in which case (\ref{TAN})
simplifies to $\lambda_1 = x_1 + y_1/x_1$.  
% \\[0,3cm]
\clearpage
\noindent 
{\bf Figure 1.} Algorithm for computing $2P + Q$  and the equation for a 
parabola through $P$, $P$, $Q$, and $-(2P+Q)$, where $P = (x_1, \, y_1)$ 
and $Q = (x_2, \, y_2)$.

\vskip .1 truein

%\begin{figure}
%Here is the pseudocode for implementing the savings. It includes the
%details for characteristic 2 and all the special cases that need to be
%considered.
%
%Let $P$ and $Q$ be points on an elliptic curve
%  
%   \[ y^2 + xy = x^3 + ax^2 + b \;\;\; \mbox{(characteristic 2)}\] or
%\[ y^2 = x^3 + ax + b \;\;\; \mbox{(field $\F_q$, $q$ odd)}\]
% 
%\noindent u\\
% \phantom{elseblah}$y_4 = \lambda_2 (x_1 - x_4) - a_1 x_4 - a_3 - y_1$;
%\\
%\phantom{else}end if\\
%\phantom{else}return $(x_4,\, y_4)$;\\
%end if
%sing affine coordinates. In the following $x(P)$ denotes
%the $x$-coordinate of a point $P$ on the elliptic curve, and $y(P)$
%denotes the $y$-coordinate.
%
\noindent if ($P = \mathbf{O}$) then \\
\phantom{else}if ($Q = \mathbf{O}$) then\\
\phantom{elsebla}$\text{parabola} = 1$; \\
\phantom{else}else \\
\phantom{elsebla}$\text{parabola} = x - x_2$; \\
\phantom{else}end if \\
\phantom{else}return $Q$;\\
else if ($Q = \mathbf{O}$) then\\
\phantom{else}if (denominator of (\ref{TAN}) is zero) then \\
\phantom{elsebla}$\text{parabola} = x - x_1$; \\
\phantom{elsebla}return $\mathbf{O}$; \\
\phantom{else}end if \\
\phantom{else}Get tangent slope $\lambda_1$ from (\ref{TAN}); \\
\phantom{else}$\text{parabola} = y - y_1 - \lambda_1 (x - x_1)$; \\
\phantom{else}$x_3 = \lambda_1 (\lambda_1 + a_1) - a_2 - 2x_1$; \\
\phantom{else}$y_3 = \lambda_1 (x_1 - x_3) - a_1 x_3 - a_3 - y_1$; \\
\phantom{else}return $(x_3,~y_3)$;\\
else\\
 \phantom{else}if ($x_1 \neq x_2$) then\\
 \phantom{elsebla}$\lambda_1 = (y_1 - y_2) / (x_1 - x_2)$; \quad
                     {\tt/* slope of line through $P$, $Q$. */}\\
 \phantom{else}else if ($y_1 \neq y_2$  OR  denominator of (\ref{TAN}) is zero) then\\
 \phantom{elsebla}$\text{parabola} = (x - x_1)^2$; \\ 
 \phantom{elsebla}return $P$; \quad
           {\tt/* $P$ and $Q$ must be negatives, so $2P + Q = P$.*/}\\
 \phantom{else}else \\
 \phantom{elsebla}Get tangent slope $\lambda_1$ from (\ref{TAN}); \\
 \phantom{else}end if\\ 
 \phantom{else}$x_3 = \lambda_1 (\lambda_1 + a_1) - a_2 - x_1 - x_2$;
\\
 \phantom{elseblah}
          {\tt/* Think $y_3 = \lambda_1 (x_1 - x_3) -a_1 x_3 - a_3 - y_1$. */}\\
 \phantom{else}if ($x_3 = x_1$) then\\
 \phantom{elsebla}$\text{parabola} = y - y_1 - \lambda_1 (x - x_1)$; \\
 \phantom{elsebla}return $\mathbf{O}$; \quad{\tt/* $P + Q$ and $P$ are negatives. */}\\
\phantom{else}end if\quad{\tt/* Think $\lambda_2 =
                        (y_1 - y_3) / (x_1 - x_3)$ */}\\
 \phantom{else}$\lambda_2 = (a_1 x_3 + a_3 + 2 y_1)/(x_1 - x_3) -
\lambda_1$; \\ 
 \phantom{else}$x_4 = \lambda_2(\lambda_2 + a_1) - a_2 - x_1 - x_3$; \\
 \phantom{else}$y_4 = \lambda_2 (x_1 - x_4) - a_1 x_4 - a_3 - y_1$; \\
 \phantom{else}$\text{parabola} = (x - x_1)
                                  (x - x_4 + (\lambda_1 + \lambda_2 + a_1) \lambda_2)
                                     - (\lambda_1 + \lambda_2 + a_1)(y - y_1)$; \\
\phantom{else}return $(x_4,\, y_4)$;\\
\noindent end if
%\caption {\bf Algorithm for computing $2P + Q$  and the equation for a 
%parabola through $P$, $P$, $Q$, and $-(2P+Q)$,
%where $P = (x_1, \, y_1)$ and $Q = (x_2, \, y_2)$.} 
%\label{FigCode}

% With \clearpage rather than \pagebreak,
% we don't get the underflow \vbox messages.
% The white space is at the bottom of the last page, not spread out.
%\end{figure}

\end{document}